\documentclass
{article}
\usepackage[T2A]{fontenc}
\usepackage[cp1251]{inputenc}
\usepackage[english,russian]{babel}
\usepackage{amsmath}
\usepackage{amsfonts,amssymb,mathrsfs,amscd}
\usepackage{verbatim}
\usepackage{eucal}
\usepackage{graphicx}
\usepackage{enumerate}
\usepackage{esint}
\usepackage{upgreek}
\let\No=\textnumero

\usepackage[hyper]{izv2e}
\JournalName{}

\numberwithin{equation}{section}

\theoremstyle{plain}
\newtheorem{maintheorem}{Основная теорема}

\newtheorem{theo}{Теорема} 

\theoremstyle{definition}
\newtheorem{proof}{Доказательство}
\newtheorem{remark}{Замечание}

\renewcommand{\leq}{\leqslant} 
\renewcommand{\geq}{\geqslant}
\newcommand{\RR}{\mathbb{R}} 
\newcommand{\CC}{\mathbb{C}} 
\newcommand{\NN}{\mathbb{N}} 

\newcommand{\DD}{\mathbb{D}} 

\DeclareMathOperator{\rad}{rad}
\DeclareMathOperator{\bal}{bal}

\DeclareMathOperator{\supp}{supp} 
\DeclareMathOperator{\type}{type} 
 
\DeclareMathOperator{\ord}{ord}
\DeclareMathOperator{\rh}{rh}
 \DeclareMathOperator{\lh}{lh}
 \DeclareMathOperator{\up}{up}
\DeclareMathOperator{\dd}{d}
\renewcommand{\Re}{\operatorname{Re}}
\renewcommand{\Im}{\operatorname{Im}}

\begin{document} 

\title{Условие Линделёфа для распределения зарядов и выметание}
\author[B.\,N.~Khabibullin]{Б.\,Н.~Хабибуллин}
\address{Башкирский государственный университет}
\email{khabib-bulat@mail.ru} 

\udk{517.538 : 517.574}

\maketitle
	
	\begin{fulltext}

\begin{abstract} 
Пусть $\nu$ --- распределение зарядов на комплексной плоскости $\mathbb C$, т.е. вещественная мера Радона на $\mathbb C$ полной вариации $|\nu|$. Распределение зарядов $\nu$ конечной верхней плотности при порядке $1$, если 
$$
\limsup_{0<r\to+\infty} \frac{1}{t}|\nu|\Bigl(\bigl\{z\in \mathbb C \bigm| |z|\leq r\bigr\}\Bigr)<+\infty.
$$ 
Распределение зарядов $\nu$ удовлетворяет условию Линделёфа рода $1$, если 
$$
\limsup_{1\leq r\to +\infty}\biggl|\int_{1\leq |z|\leq r} \frac{1}{z}\operatorname{d}\!\nu(z)\biggr|<+\infty.
$$
Эти понятия играют ключевую роль при исследованиях  целых функций экспоненциального типа и субгармонических функций  конечного типа при порядке $1$, а также при их применениях. В наших предыдущих работах была разработана техника выметания конечного 
рода $q=0,1,\dots$ распределения зарядов $\nu$ из  полуплоскости. Мы показываем, что  выметание  рода $q=1$ распределения зарядов $\nu$ из  полуплоскости  сохраняет   эту пару свойств  при некотором условии  на  $\nu$.

Библиография:   14 наименований

\end{abstract}
	

\markright{Условие Линделёфа для распределения зарядов и выметание}

\footnotetext[0]{Исследование выполнено за счёт гранта Российского научного фонда №~22-21-00026, 
\href{https://rscf.ru/project/22-21-00026/}{https://rscf.ru/project/22-21-00026/}.}


\section{Некоторые  обозначения, определения, соглашения}\label{11def} 
 Одноточечные множества $\{a\}$ часто записываем без фигурных скобок, т.е. просто как $a$. Так,   $\NN_0:=0\cup \NN=\{0,1, \dots\}$ для  множества $\mathbb N:=\{1,2, \dots\}$ {\it натуральных чисел,\/}  
$\CC_{\infty}:=\CC\cup \infty$ и $\overline \RR:=-\infty \cup \RR\cup +\infty$ --- {\it расширенные\/}
комплексная плоскость и вещественная ось с $-\infty:=\inf \RR\notin \RR$, $+\infty:=\sup \RR\notin \RR$, неравенствами $-\infty\leq x\leq +\infty$ для любого $x\in \overline \RR$
 и естественной порядковой топологией. По определению 
$\sup \varnothing:=-\infty$ и $\inf \varnothing:=+\infty$ для {\it пустого множества\/} $\varnothing$.
Символом  $0$   могут обозначаться {\it нулевые\/} функции,  меры  и пр.

Для $x\in X\subset \overline \RR$ его {\it положительную часть\/} обозначаем как $x^+:=\sup\{0,x \}$, $X^+:=\bigl\{x^+\bigm|  x\in X\bigr\}$. {\it Расширенной числовой функции\/} $f\colon S\to \overline \RR$ сопоставляем её {\it положительную часть\/} $f^+\colon s\underset{s \in S}{\longmapsto} (f(s))^+\in \overline{\RR}^+$ 
и {\it отрицательную часть\/} $f^-:=(-f)^+\colon S\to \overline{\RR}^+ $. Как обычно, пишем $f\not\equiv c$, если функция $f$ принимает хотя бы одно значение, отличное от $c$, в области её определения. 

Для $x_0\in \RR$ и расширенной числовой функции $m\colon [x_0,+\infty) \to \overline \RR$ определим  
\begin{equation}\label{senu0:a} 
\ord[m]:=\limsup_{x\to +\infty} \frac{\ln \bigl(1+m^+(x)\bigr)}{\ln x}\in
\overline \RR^+
 \end{equation} 
{\it --- порядок\/}  (роста) функции $m$ (около $+\infty$), а  для $p\in \RR^+$
\begin{equation} 
\type_p[m]:=\limsup_{x\to +\infty} \frac{m^+(x)}{x^p}\in \overline \RR^+
\label{typevf}
 \end{equation} 
{\it --- тип\/} (роста) функции $m$ {\it при порядке\/} $p$ (около $+\infty$) \cite{Boas}, \cite{Levin56}, \cite{Levin96}, 
\cite{Kiselman}, \cite[2.1]{KhaShm19}, 
а для произвольной 
 функции $u\colon \CC\to \overline \RR$ с  {\it радиальный функцией роста\/}
\begin{equation} 
{\mathrm M}_u \colon r\underset{r\in \RR^+}{\longmapsto}
\sup\bigl\{u(z)\bigm| |z|=r\bigr\}
\label{u}
\end{equation} 
по определению $\ord[u]:=\ord[{\mathrm M}_u]$ и $\type_p[u]:=\type_p[{\mathrm M}_u]$ 
--- соответственно {\it порядок\/} и 
{\it тип функции $u$ при порядке $p$} \cite{Boas}, \cite{Levin96}, \cite{Kiselman}, \cite[Замечание 2.1]{KhaShm19}.  Функции $u$ конечного типа 
$\type_1[u]\in \RR^+$ при порядке $p=1$ называем просто функциями {\it конечного типа,\/} не упоминая  и не указывая 
порядок $1$ в $\type[u]:=\type_1[u]$.

{\it Распределением масс\/} --- это  {\it положительная   мера Радона\/}  \cite{EG},  \cite[Дополнение A]{Rans}, \cite[гл.~3]{HK}, 
а {\it распределение зарядов\/} ---  разность распределений масс \cite{Landkof}.
Для распределений масс   или зарядов  {\it на\/} $\CC$ часто не указываем, где  они заданы.   
Для  субгармонической  в области из $\CC$ функции $u\not\equiv -\infty$ действие  на неё   {\it оператора Лапласа\/} ${\bigtriangleup}$  в смысле теории обобщённых функций  определяет её {\it распределение масс Рисса\/}
$\frac{1}{2\pi}{\bigtriangleup}u$  в этой области  \cite{HK}, \cite{Rans}, \cite{Az}.

 Далее $D_z(r):=\bigl\{w \in \CC \bigm| |w-z|<r\bigr\}$ и  $\overline{D}_z(r):=\bigl\{w \in \CC_{\infty} \bigm| |w-z|\leq r\bigr\}$,
а также $\partial \overline{D}_z(r):=\overline{D}_z(r)\setminus {D}_z(r)$
---  соответственно {\it открытый} и   {\it замкнутый круги,\/} а также  {\it окружность  радиуса\/ $r\in \overline \RR^+$
 с центром\/} $z\in \CC$, а   $\DD:=D_0(1)$ и  $\overline \DD:=\overline D_0(1)$, а также  
$\partial \overline \DD:=\partial \overline{D}_0(1)$ --- соответственно {\it открытый\/} и {\it замкнутый 
единичные круги,\/} а также {\it единичная окружность} в  $\CC$. 

Через $\CC_{ \rh}:=\bigl\{z\in \CC \bigm| \Re z>0\bigr\}$ и $\CC_{\overline  \rh}:=\CC_{ \rh}\cup i\RR$, а также
$\CC_{ \lh}:=-\CC_{ \rh}$ и $\CC_{\overline \lh}:=-\CC_{\overline  \rh}$ обозначаем 
соответственно {\it правые открытую} и   {\it замкнутую полуплоскости,\/} а также  
{\it левые открытую} и   {\it замкнутую  полуплоскости\/} в  $\CC$. 

Для распределения зарядов  $\nu$ на $S\subset \CC$ через $\nu^+:=\sup\{\nu,0\}$, $\nu^-:=(-\nu)^+$
и $|\nu|:=\nu^++\nu^-$ обозначаем соответственно {\it верхнюю, нижнюю\/}
и {\it полную вариации\/} распределения зарядов  $\nu$, а $\supp \nu=\supp |\nu|$ --- его {\it носитель,} но  распределение зарядов  $\nu$ {\it сосредоточено на $\nu$-измеримом подмножестве $S_0\subset S$,\/} если полная вариация $|\nu|$ дополнения $S\setminus S_0$ множества $S$ равна нулю. 

{\it Сужение\/} функции $f$ на  $S\subset \CC$ обозначаем как   $f{\lfloor}_S$. Аналогично через $\nu{\lfloor}_S$
обозначается  обозначается и   {\it сужение\/}  положительной меры Бореля или  распределения зарядов $\nu$ на  $\nu$-измеримое  $S\subset \CC$.
При $r\in \overline \RR^+$ для таких   $\nu$  через 
\begin{equation}\label{df:nup} 
\nu_z^{\rad} (r):=\nu \bigl(\,\overline D_z(r)\bigr),\quad \nu^{\rad}(r):=\nu_0^{\rad}(r)
=\nu\bigl(r\overline \DD\bigr)
\end{equation} 
обозначаем  {\it радиальные непрерывные справа считающие  функции  распределения зарядов  $\nu$ с
центрами\/}  соответственно {\it в точке\/ $z\in \CC$} и {\it в нуле.\/}

{\it Верхняя плотность распределения зарядов\/ $\nu$ при порядке\/} $p\in \RR^+$  равна 
\begin{equation}\label{typenu}
\type_p[\nu]:=\type_p\bigl[|\nu|\bigr]
\overset{\eqref{typevf}}{:=} \limsup_{0<r\to +\infty} \frac{|\nu|(r\overline \DD)}{r^p}
\overset{\eqref{df:nup}}{=} \limsup_{0<r\to +\infty} \frac{|\nu|^{\rad}(r)}{r^p}
\in \overline \RR^+,
\end{equation} 
и при $p=1$ упоминание о порядке  опускаем.   В частности, {\it  распределение зарядов\/} $\nu$
{\it конечной верхней плотности,\/} если $\type[\nu]:=\type_1[\nu]<+\infty$. {\it Порядок распределения зарядов\/}
$\nu$ определяется как $\ord[\nu]\overset{\eqref{senu0:a}}{:=}\ord\bigl[|\nu|^{\rad}\bigr]$
через \eqref{df:nup}.

\section{Выметание распределения зарядов  из  правой полуплоскости}\label{BalCrh}

Для распределения заряда $\nu$ используем \cite[формула (1.9)]{KhaShm19}   его {\it функцию  распределения на\/}  $\RR$, 
обозначаемую как  $\nu_{\RR}\colon \RR\to  \RR$ и определённую  равенствами 
\begin{equation}\label{nuR} 
\nu_{\RR}(x_2)-\nu_{\RR}(x_1):=\nu\bigl((x_1,x_2]\bigr), \quad
-\infty <x_1<x_2<+\infty, 
\end{equation}
и   {\it функцию распределения $\nu_{i\RR}\colon \RR\to \RR$    на\/} $i\RR$, определённую равенствами
\begin{equation}\label{nuiR} 
\nu_{i\RR}(y_2)-\nu_{i\RR}(y_1):=\nu\bigl(i(y_1,y_2]\bigr), \quad -\infty <y_1<y_2<+\infty. 
\end{equation}
Поскольку эти функции распределения определены лишь с точностью до аддитивной постоянной, 
при необходимости используем их {\it нормировки в нуле}
\begin{equation}\label{nuo}
\nu_{\RR}(0):=0, \quad \nu_{i\RR}(0):=0.
\end{equation}
По построению \eqref{nuR} и \eqref{nuiR}    функции  $\nu_{\RR}$ и $\nu_{i\RR}$  локально ограниченной вариации  на $\RR$, а в случае {\it распределения масс\/} $\nu$ обе эти функции возрастающие. 
Обратно, любая функция локально ограниченной вариации на $\RR$ или $i\RR$ однозначно определяет распределение зарядов с носителем соответственно на  $\RR$ или $i\RR$. 

Мы напоминаем и адаптируем основные понятия и утверждения из \cite{KhaShm19} и \cite{KhaShmAbd20},
 а также частично из \cite{Kha91} и \cite{Kha91AA} о  выметании конечного рода $q\in \NN_0$ распределений зарядов, но пока применительно только к правой полуплоскости $\CC_{\rh}$ в случаях  $q:=0$ и $q:=1$. В 
\cite{KhaShm19} и \cite{KhaShmAbd20} в основном рассматривается {\it верхняя полуплоскость\/} $\CC^{\up}:=i\CC_{ \rh}$, что переносится на $\CC_{\rh}$ поворотом на прямой угол.

{\it Характеристическую функцию множества\/}  $S$ обозначаем через 
\begin{equation}\label{SdrS}
\boldsymbol{1}_S\colon z\underset{z\in \mathbb C}{\longmapsto} \begin{cases}
1&\text{ если $z\in S$},\\
0&\text{ если $z\notin S$}.
\end{cases}
\end{equation}
 
{\it Гармоническая  мера  для\/ $\CC_{\rh}$ в точке\/ $z\in \CC_{\rh}$} на интервалах $i(y_1,y_2]\subset i\overline\RR$
\begin{equation}\label{omega}
\omega_{\rh} \bigl(z,i(y_1,y_2]\bigr){\overset{\text{\cite[3.1]{KhaShm19}}}{:=}}\omega_{\CC_{\rh}}(z,i(y_1,y_2])
\underset{z\in \CC_{\rh}}{:=}\frac1{\pi}
\int_{y_1}^{y_2}\Re \frac{1}{z-iy} \dd y 
\end{equation} 
равна делённому на $\pi$ углу, под которым виден интервал $i(y_1,y_2]$ из точки $z\in \CC_{\rh}$ \cite[(3.1)]{Kha91AA}, \cite[1.2.1, 3.1]{KhaShm19},  а  в точках  мнимой оси $i\RR$ определяется как 
\begin{equation}\label{oiR}
\omega_{\rh} \bigl(iy,i(y_1,y_2]\bigr):=\boldsymbol{1}_{(y_1,y_2]}(y)
\quad\text{при $y\in  \RR$.}
\end{equation}
Для  распределения зарядов  $\nu $ при {\it классическом условии Бляшке\/} для  $\CC_{\rh}$
\begin{equation}\label{Blcl}
 \int_{\CC_{\rh}\setminus  \DD} \Re \frac{1}{z}\dd |\nu| (z)<+\infty
\end{equation}
определено \cite[следствие 4.1, теорема 4]{KhaShm19}  
его классическое выметание из $\CC_{\rh}$ на  $ \CC_{\overline \lh}$ с носителем на  $  \CC_{\overline \lh}$, которое  в более широких рамках  \cite[определение 3.1]{KhaShmAbd20} 
представляет собой {\it выметание рода\/ $0$,\/} обозначавшееся в \cite{KhaShmAbd20} как
$\nu^{\bal[0]}_{\CC_{\overline \lh}}$. 
Здесь используется чуть более компактная  запись  
 $\nu^{\bal^0_{\rh}}:=\nu^{\bal[0]}_{\CC_{\overline \lh}}$. По определению {\it распределение зарядов\/}  $\nu^{\bal^0_{\rh}}$ --- это {\it сумма сужения\/}  $\nu{\lfloor}_{\CC_{\lh}}$ на $\CC_{\lh}$  с {\it распределением зарядов на $i\RR$,} определяемым в обозначениях \eqref{nuiR}  {\it функцией распределения\/}
\begin{equation}\label{mubal}
\nu^{\bal^0_{\rh}}_{i\RR}(y_2)-\nu^{\bal^0_{\rh}}_{i\RR}(y_1)\overset{\eqref{omega},\eqref{oiR}}{:=}
\int\limits_{\CC_{\overline \rh}} \omega_{\rh}\bigl(z, i(y_1,y_2]\bigr) \dd \nu(z)
\end{equation}
с нормировкой вида \eqref{nuo} при необходимости. 
Классическое выметание  рода $0$ не увеличивает полную  меру полной вариации  распределения зарядов, 
поскольку гармоническая мера \eqref{omega} вероятностная и 
\begin{equation}\label{omega1}
\bigl|\nu^{\bal^0_{\rh}}\bigr|(S)\underset{S\subset \CC}{\overset{\eqref{mubal}}{\leq}} |\nu|(S).
\end{equation}  

В \cite[определение 2.1]{KhaShmAbd20} вводилось понятие гармонического   заряда  рода\/ $1$ для верхней полуплоскости \/ $\CC^{\up}$ в точке\/ $z\in \CC^{\up}$,  обозначавшегося  в \cite[формула (2.1)]{KhaShmAbd20} через  $\Omega^{[1]}_{\CC^{\up}}$. Здесь используем 
поворот на прямой угол с переходом от $\CC^{\up}$ к $\CC_{\rh}$ и определим 
 {\it гармонический    заряд рода\/ $1$ для правой  полуплоскости \/ $\CC_{\rh}$} 
 как функцию $\Omega_{\rh}$ ограниченных  интервалов $i(y_1,y_2]\subset i\RR$ по правилу 
\begin{equation}\label{Ocr}
\Omega_{\rh}\bigl(z,i(y_1,y_2]\bigr)
\overset{\eqref{omega},\eqref{oiR}}{:=}\omega_{\rh}\bigl(z,i(y_1,y_2]\bigr)-\frac{y_2-y_1}{\pi}\Re\frac{1}{z}
\quad\text{в $z\in \CC_{\overline \rh}\setminus 0$}.
\end{equation}
Для  распределения зарядов  $\nu$ в  \cite[определение 3.1, теорема 1, замечание 3.3]{KhaShmAbd20} определялось 
{\it выметание\/ $\nu_{\CC_{\overline \lh}}^{\bal[1]}$ рода\/ $1$ распределения зарядов\/  $\nu$ из\/  $\CC_{\rh}$ на \/} $ \CC_{\overline \lh}$
  при $0\notin \supp \nu$. Здесь используется  более компактная запись  для такого выметания $\nu^{\bal^1_{\rh}}:=\nu^{\bal[1]}_{ \CC_{\overline \lh}}$.  По определению {распределение зарядов\/} $\nu^{\bal^1_{\rh}}$ --- это 
{\it сумма  сужения\/} $\nu{\lfloor}_{\CC_{\lh}}$ с {\it  распределением зарядов на\/} $i\RR$, определяемым в обозначениях \eqref{nuiR}  {\it функцией распределения} 
\begin{equation}\label{df:nurh}
\nu^{\bal^1_{\rh}}_{i\RR}(y_2)-\nu^{\bal^1_{\rh}}_{i\RR}(y_1)\overset{\eqref{Ocr}}{=}
\int_{\CC_{\rh}} \Omega_{\rh} \bigl(z, i(y_1,y_2]\bigr)\dd \nu (z)
\end{equation}
с нормировкой вида \eqref{nuo} при необходимости. 
\begin{remark}\label{remAA} Выметание распределения зарядов  рода $1$ из $\CC_{\rh}$ на  $ \CC_{\overline \lh}$ --- это часть глобального, или двухстороннего, выметания распределения зарядов  из $\CC\setminus i\RR$ на $i\RR$, рассмотренного  в \cite[\S~3]{Kha91AA} и сыгравшего там одну из ключевых ролей. Двустороннее выметание на мнимую ось можно рассматривать как последовательное выметание рода $1$ распределения зарядов  сначала из $\CC_{\rh}$ на $ \CC_{\overline \lh}$, а затем зеркально симметричной относительно $i\RR$ процедуры выметания рода $1$ получившегося распределения зарядов  из  $\CC_{\lh}$ на  $\CC_{\overline \rh}$. 
\end{remark}

Ограничение $0\notin \supp \nu$ для выметания рода $1$ легко преодолевается, если скомбинировать выметание рода $0$ части $\nu$  около нуля с выметанием рода $1$ для оставшейся части $\nu$. Для этого определяем {\it комбинированное выметание рода\/ $01$ 
распределения зарядов\/ $\nu$ из\/ $\CC_{\rh}$ на\/ $\CC_{\overline \lh}$} \cite[замечание 3.3, (3.43), (4.1)]{KhaShmAbd20}
\begin{equation}\label{bal01}
\nu^{\bal_{\rh}^{01}}:=\bigl(\nu{\lfloor}_{r_0\DD}\bigr)^{\bal_{\rh}^{0}}+\bigl(\nu{\lfloor}_{\CC\setminus r_0\DD}\bigr)^{\bal_{\rh}^{1}}
\end{equation}
при каком-нибудь фиксированном радиусе  $r_0\in \RR^+\setminus 0$. 
\begin{remark}\label{rem0nu}
Круг $r_0\DD$  в правой части \eqref{bal01}  можно заменить на любое ограниченное борелевское множество, содержащее полукруг $r_0\DD\setminus \CC_{\overline \lh}$, или даже на пустое множество.  Другими словами, можно  обойтись совсем без 
выметания рода $0$ и положить $\nu^{\bal_{\rh}^{01}}:=\nu^{\bal_{\rh}^{1}}$
если $0\notin \supp \nu$ или, более общ\'о,
\begin{equation*}
\int_{r_0\DD} \Re^+\frac{1}{z} \dd |\nu|(z)<+\infty, 
\end{equation*}
что, очевидно, выполнено, если для некоторого $r_0\in \RR^+\setminus 0$ имеем 
\begin{equation}\label{nur0+}
|\nu|\bigl(r_0\DD\cap \CC_{ \rh}\bigr)=0.
\end{equation}
\end{remark}
Теперь вопрос существования выметания  $\nu^{\bal_{\rh}^{01}}$  упирается лишь в поведение 
распределения зарядов $\nu$ около бесконечности.

Распределение зарядов $\nu$ принадлежит  {\it классу сходимости при порядке\/} роста  $p\in \NN_0$, если
\cite[определение 4.1]{HK}, \cite[\S~2, 2.1, (2.3)]{KhaShm19} 
\begin{equation}\label{sufc}
\int_1^{+\infty}\frac{|\nu|^{\rad}(t)}{t^{p+1}}\dd t<+\infty. 
\end{equation}
Мы  используем  различные виды условий Линделёфа. 
Распределение зарядов $\nu$ удовлетворяет  {\it $\RR$-условию Линделёфа\/}  (рода $1$), если  
\begin{equation}
\sup_{r\geq 1} \biggl| \int_{1<|z|\leq r}\Re\frac{1}{z}\dd \nu(z)\biggr|<+\infty, 
\label{con:LpZR}
\end{equation}
удовлетворяет {\it $i\RR$-условию Линделёфа\/} (рода $1$), если
\begin{equation}
\sup_{r\geq 1}\biggl| \int_{1<|z|\leq r}\Im\frac{1}{z}\dd \nu(z)
\biggr|<+\infty, 
\label{con:LpZiR}
\end{equation}
и удовлетворяет  {\it условию Линделёфа\/} (рода $1$), если
\begin{equation}
\sup_{r\geq 1}\biggl|\int_{1<|z|\leq r}\frac{1}{z}\dd \nu(z)
\biggr|<+\infty. 
\label{con:LpZ}
\end{equation}

Ключевая роль условий Линделёфа отражает следующая классическая  
\begin{theo}[{Вейерштрасса\,--\,Адамара\,--\,Линделёфа\,--\,Брело (\cite[3, Теорема  12]{Arsove53p},  \cite[4.1, 4.2]{HK}, \cite[2.9.3]{Az}, \cite[6.1]{KhaShm19})}]\label{pr:rep}
Если   $u\not\equiv -\infty$ --- субгармоническая функция конечного типа, то её   распределение масс Рисса  $\frac{1}{2\pi}{\bigtriangleup}u$
конечной верхней плотности и удовлетворяет условию Линделёфа \eqref{con:LpZ}. 

Обратно, если  распределение масс    $\nu$ конечной верхней плотности, то существует  субгармоническая функция $u_{\nu}$ с  
$\frac{1}{2\pi}{\bigtriangleup}u_{\nu}=\nu$  порядка $\ord[u_{\nu}]
\leq 1$, которая при выполнении условия Линделёфа \eqref{con:LpZ} для\/  $\nu$   будет уже  функцией  конечного типа. 
При этом  любая  субгармоническая   функция $u$ с  $\frac{1}{2\pi}{\bigtriangleup}u=\nu$
представляется в виде суммы $u=u_{\nu}+H$, где $H$ --- гармоническая функция на\/ $\CC$, которая при условии\/ $\type_2[u]=0$ является  гармоническим многочленом степени\/  $\deg H\leq 1$, а функция $u$ становится функцией порядка\/  $\ord[u]\leq 1$.
\end{theo}

В свете теоремы Вейерштрасса\,--\,Адамара\,--\,Линделёфа\,--\,Брело важна 

\begin{maintheorem}\label{theoB1} Пусть $\nu$  ---  распределение зарядов, для которого сужение $\nu{\lfloor}_{\CC_{ \rh}}$ при\-н\-а\-д\-л\-ежит   классу сходимости при порядке $p\overset{\eqref{sufc}}{=}2$.  
Тогда существует выметание $\nu^{\bal^{01}_{\rh}}$ рода $01$ из $\CC_{\rh}$ на $ \CC_{\overline \lh}$. В частности, если $\ord[\nu]<2$, то  $\nu$ из  класса сходимости при порядке $<2$ и   $\ord[\nu^{\bal^{01}_{\rh}}]\leq \ord[\nu]$.  

Если  $\nu$ конечной верхней плотности и удовлетворяет условию 
\begin{equation}\label{cB2l}
\sup_{r\geq 1}\biggl|\int_{1 < | z|\leq r} \Re^+ \frac{1}{ z} \dd \nu(z)\biggr|<+\infty,
\end{equation}
то $\nu^{\bal^{01}_{\rh}}$ ---  распределение зарядов  конечной верхней плотности  c разностью $\nu-\nu^{\bal^{01}_{\rh}}$,  удовлетворяющей условию Линделёфа. При этом  если   носитель $\supp \nu$ не пересекается с замкнутым углом раствора строго больше, чем  $\pi$, содержащим  $i\RR$,  с вершиной в нуле и биссектрисой $\RR$, т.е.
\begin{equation}\label{eA1}
 \bigl\{z\in \CC \bigm| \Re z\leq a |z|\bigr\}  \bigcap \supp \nu=\varnothing 
\quad\text{для  некоторого $a \in (0,1)$},
\end{equation}
то в обозначении \eqref{df:nup} для радиальных считающих функций $|\nu^{\bal^{01}_{\rh}}\bigr|_{iy}^{\rad}$ полной вариации $|\nu^{\bal^{01}_{\rh}}|$ распределения зарядов  $\nu^{\bal^{01}_{\rh}}$ с центрами  $iy\in i\RR$  имеем
\begin{equation}\label{trnuair}
\sup_{y\in \RR}\sup_{t\in (0,1]}\frac{|\nu^{\bal^{01}_{\rh}}\bigr|_{iy}^{\rad}(t)}{t}<+\infty.
\end{equation}

\end{maintheorem}

\begin{proof} 
Для классического выметания $\bigl(\nu{\lfloor}_{r_0\overline \DD}\bigr)^{\bal_{\rh}^{0}}$
рода $0$ из правой части \eqref{bal01}  имеем 
$$
\Bigl|\bigl(\nu{\lfloor}_{r_0\DD}\bigr)^{\bal_{\rh}^{0}}\Bigr|(\CC)\overset{\eqref{omega1}}{\leq} |\nu|\bigl(r_0\DD\bigr).
$$ 
Отсюда, очевидно,
$\bigl(\nu{\lfloor}_{r_0\DD}\bigr)^{\bal_{\rh}^{0}}$ конечной верхней плотности и удовлетворяет условию Линделёфа,
 поэтому далее можем ограничиться распределениями зарядов $\nu$ с носителем 
\begin{equation}\label{2as}
\supp \nu \subset \CC\setminus 2d\overline \DD \quad\text{для  некоторого $d\in (0,1]$}
\end{equation}
и рассматривать  
 только выметание $\nu^{\bal^1_{\rh}}$ рода $1$. Исходя из этого, 
первые утверждения об условии существования  выметания $\nu^{\bal^{01}_{\rh}}$ и его порядке 
  --- очень частные случаи \cite[теорема 1]{KhaShmAbd20}, а  
конечность верхней плотности  выметания $\nu^{\bal^1_{\rh}}$ при дополнительном 
условии \eqref{cB2l}  --- частный случай  \cite[теорема 3, п.~4]{KhaShmAbd20} или же  \cite[теорема 3.1]{Kha91AA}, если учесть замечание 
\ref{remAA}.

Докажем при условии \eqref{cB2l} выполнение условия Линделёфа для разности $\nu -\nu^{\bal^1_{\rh}}$, что будет развитием и усилением  \cite[теорема 3.2]{Kha91AA}, где это доказывалось при дополнительном требовании, что и само   
распределение зарядов  $\nu$  удовлетворяет условию  Линделёфа. 

Не умаляя общности можно считать, что  {\it распределение зарядов  $\nu$ сосредоточено в\/ $\CC_{\rh}$,} 
поскольку сужение $\nu{\lfloor}_{\CC_{\overline \lh}}$ при выметании из $\CC_{ \rh}$ не меняется и 
$$
\nu-\nu^{\bal_{\rh}^1}=\nu{\lfloor}_{\CC_{ \rh}}-\bigl(\nu{\lfloor}_{\CC_{ \rh}}\bigr)^{\bal_{\rh}^1}.
$$
Когда $\nu$ сосредоточено в $\CC_{ \rh}$, имеем  $\supp \nu^{\bal^1_{\rh}}\subset  i\RR$ и из  \eqref{cB2l} следует
\begin{equation}\label{renu}
\biggl|\int_{(r\overline \DD)\setminus d\overline \DD}\Re \frac{1}{z}\dd \bigl(\nu-\nu^{\bal^1_{\rh}}\bigr)(z)\biggr|
\overset{\eqref{cB2l}}{\underset{r\to+\infty}{=}}O(1).
\end{equation}
Таким образом, $\RR$-условие Линделёфа  для разности зарядов $\nu-\nu^{\bal^1_{\rh}}$ выполнено и достаточно установить 
$i\RR$-условие Линделёфа  для $\nu-\nu^{\bal^1_{\rh}}$.

Ввиду \eqref{2as} существует $C\in \RR^+$, для которого  
\begin{equation}\label{Cp}
|\nu|^{\rad}(t)+\bigl|\nu^{\bal^1_{\rh}}\bigr|^{\rad}(t)\leq Ct \quad\text{при всех $t\in \RR^+$}.
\end{equation}
{\it При фиксированном\/} $r>2$  представим  распределение зарядов  $\nu$ в виде суммы двух распределений зарядов 
\begin{equation}\label{nu+}
\nu :=\nu_{2r}+\nu_{\infty}, \quad \nu_{2r}:= \nu{\lfloor}_{2r\overline \DD}, \quad
\nu_{\infty}:=\nu -\nu_{2r}=\nu{\lfloor}_{\CC\setminus 2r\overline \DD}.
\end{equation}
По определению \eqref{df:nurh} через  функцию распределения \eqref{nuiR}  имеем равенства
\begin{multline}\label{sumnu}
-\int _{(r\overline \DD)\setminus d\overline \DD}
\Im \frac{1}{z} \dd \nu^{\bal^1_{\rh}}(z)
\overset{\eqref{df:nurh}}{=}\int_{d<|y|\leq r} \frac{1}{y}\dd \bigl(\nu^{\bal^1_{\rh}}\bigr)_{i\RR}(y)
\\\overset{\eqref{nu+}}{=}\int_{d<|y|\leq r} \frac{1}{y}\dd \bigl(\nu_{2r}^{\bal^0_{\rh}}\bigr)_{i\RR}(y)
+ 
\Bigl(-\frac{1}{\pi}\int\Re\frac{1}{z}\dd \nu_{2r}(z)\Bigr) \int_{1<|y|\leq r} \frac{\dd y}{y}
\\
+\int_{d<|y|\leq r} \frac{1}{y}\dd \bigl(\nu_{\infty}^{\bal^1_{\rh}}\bigr)_{i\RR}(y)\\
=
\int_{|y|>d} \frac{1}{y}\dd \bigl(\nu_{2r}^{\bal^0_{\rh}}\bigr)_{i\RR}(y)
+\int_{|y|>r} -\frac{1}{y}\dd \bigl(\nu_{2r}^{\bal^0_{\rh}}\bigr)_{i\RR}(y)
\\+ \int_{d<|y|\leq r} \frac{1}{y}\dd \bigl(\nu_{\infty}^{\bal^1_{\rh}}\bigr)_{i\RR}(y).
\end{multline}
Для второго интеграла в правой части \eqref{sumnu}  получаем неравенства  
\begin{equation}\label{y>r}
\Bigl|\int_{|y|>r} -\frac{1}{y}\dd \bigl(\nu_{2r}^{\bal^0_{\rh}}\bigr)_{i\RR}(y)\Bigr|
\overset{\eqref{nu+}}{\leq} \frac{1}{r}\bigl|\nu_{2r}^{\bal^0_{\rh}}\bigr|^{\rad}(2r) 
\overset{\eqref{omega1}}{\leq} \frac{1}{r}|\nu_{2r}|^{\rad}(2r) 
\overset{\eqref{Cp}}{\leq} 2C.
\end{equation}
Для последнего интеграла в правой части \eqref{sumnu} имеем 
\begin{multline}\label{kpr}
\Bigl|\int_{d<|y|\leq r} \frac{1}{y}\dd (\nu_{\infty}^{\bal^1_{\rh}})_{i\RR}(y)\Bigr|\leq 
\int_d^r \frac{1}{t}\dd \bigl|\nu_{\infty}^{\bal^1_{\rh}}\bigr|^{\rad}(t)
\\
\leq \frac{\bigl|\nu_{\infty}^{\bal^1_{\rh}}\bigr|^{\rad}(r)}{r}+
\int_d^r\frac{\bigl|\nu^{\bal^1_{\rh}}\bigr|^{\rad}(t)}{t^2} \dd t
\overset{\eqref{Cp}}{\leq} C+\int_d^r\frac{\bigl|\nu^{\bal^1_{\rh}}\bigr|^{\rad}(t)}{t^2} \dd t.
\end{multline}
Для гармонического заряда $\Omega_{\rh}$ из \eqref{df:nurh} следует  \cite[лемма 3.1]{Kha91AA}
\begin{equation}\label{Om2r}
\Omega_{\rh}\bigl(z,i[-t,t]\bigr)\leq 2\frac{t^2}{|z|^2}\quad \text{при $2t\leq |z|$}.
\end{equation}
Используя это неравенство, ввиду  $\supp \nu_{\infty}\overset{\eqref{nu+}}{\subset} \CC\setminus 2r\DD$  получаем 
\begin{multline*}
\bigl|\nu_{\infty}^{\bal^1_{\rh}}\bigr|^{\rad}(t)\overset{\eqref{df:nurh}}{\leq} 
\int_{|z|\geq 2r} \bigl|\Omega_{\rh} \bigl(z, i[-t,t]\bigr)\bigr|\dd |\nu_{\infty}| (z)
\overset{\eqref{Om2r}}{\leq} 2t^2\int_{2r}^{\infty}\frac{\dd |\nu_{\infty}|^{\rad}(s)}{s^2} \\
\leq 4t^2\int_{2r}^{\infty} \frac{|\nu|^{\rad}(s)}{s^3}\dd s
\overset{\eqref{Cp}}{\leq} 4t^2\int_{2r}^{\infty} \frac{Cs}{s^3}\dd s=
2C\frac{t^2}{r} \quad\text{ при $t\leq r$}, 
\end{multline*}
что при применении к интегралу в правой части \eqref{kpr} даёт 
\begin{equation*}\label{}
\Bigl|\int_{d<|y|\leq r} \frac{1}{y}\dd (\nu_{\infty}^{\bal^1_{\rh}})_{i\RR}(y)\Bigr|
\leq C+\int_d^r\frac{2C{t^2}/{r}}{t^2} \dd t\leq 3C.
\end{equation*}
Применяя эту оценку вместе с \eqref{y>r} к \eqref{sumnu} получаем 
\begin{equation}\label{I11}
\biggl|\int_{|y|>d} \frac{1}{y}\dd \bigl(\nu_{2r}^{\bal^0_{\rh}}\bigr)_{i\RR}(y)+\int _{(r\overline \DD)\setminus d\overline \DD}
\Im \frac{1}{z} \dd \nu^{\bal^1_{\rh}}(z)\biggr|\leq 5C \quad\text{при всех $r>2$.}
\end{equation}
Для распределения зарядов  $\nu_{2r}$, сосредоточенного  в 
$\CC_{\rh}$,   с компактным носителем в $\CC_{\overline \rh}$ первый интеграл под модулем в левой части  
\eqref{I11} 
\begin{equation}\label{I12}
\int_{|y|>d} \frac{1}{y}\dd \bigl(\nu_{2r}^{\bal^0_{\rh}}\bigr)_{i\RR}(y)
\end{equation}
равен (см. \cite[теорема 1.2]{Kha91}, \cite[теорема 6]{KhaShm19},  \cite[(3.14)]{Kha91AA}) 
интегралу по заряду $\nu_{2r}$ от интеграла Пуассона в $\CC_{ \rh}$ от функции на $i\RR$, заданной как 
\begin{equation}\label{I2}
iy\underset{y\in \RR}{\longmapsto}
\begin{cases}
1/y&\text{ при $|y|>d$},\\
0 &\text{ при $|y|\leq d$},
\end{cases} 
\end{equation}
или от гармонического продолжения функции \eqref{I2} в $\CC_{ \rh}$. Это продолжение выписывается в явном виде 
\cite[(3.15)]{Kha91AA} и в точках $z\in \CC_{ \rh}$  принимает значение
\begin{equation*}
\biggl(\frac{1}{\pi}\Re \frac{1}{z} \ln\Bigl|\frac{z-id}{z+id}\Bigr|+\Im \frac{1}{z}\,\omega_{\rh} \bigl(z,   i[-d,d]\bigr)
\biggr) -\Im \frac{1}{z}.
\end{equation*}
Отсюда интеграл \eqref{I12} равен 
\begin{multline*}
\int_{\CC}\biggl(\frac{1}{\pi}\Re \frac{1}{z} \ln\Bigl|\frac{z-id}{z+id}\Bigr|+\Im \frac{1}{z}\,\omega_{\rh} \bigl(z,   i[-d,d]\bigr)
\biggr) \dd \nu_{2r}(z)-\int_{\CC}\Im \frac{1}{z}\dd \nu_{2r}(z)
\\\overset{\eqref{nu+},\eqref{2as}}{=}
\int_{2d<|z|\leq 2r}\biggl(\frac{1}{\pi}\Re \frac{1}{z} \ln\Bigl|\frac{z-id}{z+id}\Bigr|+\Im \frac{1}{z}\,\omega_{\rh} \bigl(z,   i[-d,d]\bigr)\biggr) \dd \nu(z)\\
-\int_{r<|z|\leq 2r}\Im \frac{1}{z}\dd \nu(z)-\int_{r\DD\setminus d\overline \DD}\Im \frac{1}{z}\dd \nu(z), 
\end{multline*}
а подстановка правой части в \eqref{I11} даёт  неравенства 
\begin{multline}\label{estt}
\biggl|-\int_{(r\overline \DD)\setminus d\overline \DD}\Im \frac{1}{z}\dd \nu(z)+\int _{(r\overline \DD)\setminus d\overline \DD}
\Im \frac{1}{z} \dd \nu^{\bal^1_{\rh}}(z)\biggr|
\\
\leq 5C+\biggl|-\int_{r<|z|\leq 2r}\Im \frac{1}{z}\dd \nu(z)\biggr|
\\+
\biggl|\int_{2d<|z|\leq 2r}\biggl(\frac{1}{\pi}\Re \frac{1}{z} \ln\Bigl|\frac{z-id}{z+id}\Bigr|+\Im \frac{1}{z}\,\omega_{\rh} \bigl(z,   i[-d,d]\bigr)\biggr) \dd \nu(z)\biggr|\\
\leq 5C+\frac{1}{r}|\nu|^{\rad}(2r)\\
+\int _{2d<|z|\leq 2r} \biggl(\frac{1}{\pi}\Bigl|\ln\Bigl|\frac{z-id}{z+id}\Bigr|\Bigr|+
\omega_{\rh} \bigl(z,    i[-d,d]\bigr)\biggr)\frac{\dd |\nu|(z)}{|z|}
\quad\text{при всех $r>2$.}
\end{multline}
  Поскольку при   $|z|\geq 2d$ имеем $|z\pm id|\geq |z|/2$, из неравенств 
\begin{equation*}
- \ln \Bigl(1+\frac{2d}{|z-id|}\Bigr)\leq \ln\Bigl|\frac{z-id}{z+id}\Bigr|\leq \ln \Bigl(1+\frac{2d}{|z+id|}\Bigr)
\end{equation*}
получаем 
\begin{equation}\label{es1}
\Bigl|\ln\Bigl|\frac{z-id}{z+id}\Bigr|\Bigr|\leq \ln \Bigl(1+\frac{4d}{|z|}\Bigr)\leq \frac{4d}{|z|}\leq \frac{4}{|z|}
\quad\text{при всех $|z|\geq 2d$.} 
\end{equation}
При   $|z|\geq 2d$ имеем  $|z|- d\geq |z|/2$ и  из  \cite[предложение 3.1, (3.10)]{KhaShm19} следует   
\begin{equation}\label{es2}
\omega_{\rh} \bigl(z,  i[-d,d]\bigr)\leq 
\frac{1}{\pi}\frac{2d\Re z}{|z|^2-d^2}\leq \frac{1}{\pi}\frac{4d\Re z}{|z|^2}\leq 
\frac{1}{\pi}\frac{4}{|z|}\quad\text{при всех $|z|\geq 2d$.} 
\end{equation}
Используя \eqref{es1} и \eqref{es2},  из  неравенств   \eqref{estt} получаем 
\begin{multline*}
\biggl|\int_{(r\overline \DD)\setminus d\overline \DD}\Im \frac{1}{z}\dd \bigl(\nu- \nu^{\bal^1_{\rh}}\bigr)(z)\biggr|
\\
\overset{\eqref{Cp}}{\leq} 7C +\int _{2d<|z|\leq 2r} \biggl(\frac{1}{\pi}\Bigl|\ln\Bigl|\frac{z-id}{z+id}\Bigr|\Bigr|+
\omega_{\rh} \bigl(z,  i[-d,d]\bigr)\biggr)\frac{\dd \nu(z)}{|z|}\\
\overset{ \eqref{es1},\eqref{es2}}{\leq} 
7C +\frac{8}{\pi} \int _{2d<|z|\leq 2r} \frac{\dd \nu(z)}{|z|^2}
\leq 7C +3 \int _{2d}^{+\infty} \frac{\dd |\nu|^{\rad}(t)}{t^2}
\\
\overset{\eqref{Cp}}{\leq} 7C +6 \int _{2d}^{+\infty} \frac{|\nu|^{\rad}(t)}{t^3}\dd t
\overset{\eqref{Cp}}{\leq}
 7C +6 C\int _{2d}^{+\infty} \frac{\dd t}{t^2}= 7C +3\frac{C}{d}
\text{ при всех $r\geq 2$.} 
\end{multline*}
Таким образом,  разность зарядов $\nu-\nu^{\bal^1_{\rh}}$ удовлетворяет условию Линделёфа.

Соотношение \eqref{trnuair} в условиях \eqref{cB2l} и \eqref{eA1}  --- частный случай сочетания 
\cite[следствие 4.2(ii)]{KhaShm19} для выметания рода $0$ и \cite[следствие 3.1, п.~(ii), (3.24)]{KhaShmAbd20} для выметания рода $1$, что неявно отражено и  в  \cite[теорема 3.3, (3.18)]{Kha91AA}.
\end{proof}

\end{fulltext}

\end{document}